\documentclass{tran-l}
 \newtheorem{thm}{Theorem}[subsection]
 \newtheorem{cor}[thm]{Corollary}
 
 \newtheorem{prop}[thm]{Proposition}
 \theoremstyle{definition}
 \newtheorem{defn}[thm]{Definition}
 \theoremstyle{remark}
 \newtheorem{rem}[thm]{Remark}
 \numberwithin{equation}{subsection}
 
\usepackage{graphicx}
\begin{document}

\title[Chen-Ruan Cohomology ]
 {The Chen-Ruan Cohomology of Almost Contact  Orbifolds }

\author{ Fan Ding, Yunfeng Jiang  and Jianzhong Pan}

\address{School of Mathematics, Peking University, Beijing 100871, P.R.China}

\email{dingfan@math.pku.edu.cn}

\address{Department of Mathematics, University of British Columbia, 1984 Mathematics Rd, Vancouver, BC, V6T 1Z2, Canada }

\email{jiangyf@math.ubc.ca}
\address{Institute of Mathematics, Academy of Mathematics and System
Sciences, Chinese Academy of Sciences, Beijing 100080, P.R.China}
\email{pjz@mail.amss.ac.cn}

\thanks{}

\thanks{}



\date{}

\dedicatory{}

\commby{}


\begin{abstract}
Comparing to the Chen-Ruan cohomology theory for the almost
complex orbifolds, we study the orbifold cohomology theory for
almost contact orbifolds. We define the Chen-Ruan cohomology group
of any almost contact orbifold. Using the methods for almost
complex orbifolds (see [2]), we define the obstruction bundle for
any 3-multisector of the almost contact orbifolds and the
Chen-Ruan cup product for the Chen-Ruan cohomology. We also prove
that under this cup product the direct sum of all dimensional
orbifold cohomology groups constitutes a cohomological ring.
Finally we calculate two examples.
\end{abstract}

\maketitle

\section{Introduction}
Motivated by the orbifold string theory on the global quotient
orbifolds studied by Dixon, Harvey, Vafa and Witten (see [4],
[5]), and also Zaslow [14], Chen and Ruan [2] developed a very interesting
cohomology theory for almost complex orbifolds which is now called
Chen-Ruan cohomology. The remarkable aspect of this theory,
besides the generalization of non global quotient orbifolds
studied by physicists, is the existence of a ring structure
obtained from Chen-Ruan's orbifold quantum cohomology construction
(see [3]) by restricting to the class called  ghost maps, the
same as the ordinary cup product may be obtained by quantum cup
product. The appearance of this interesting theory attracted many
authors to calculate the Chen-Ruan cohomological rings of concrete
examples, see [6],[12],[10],[7]. But we would like to point out
that the definition of the Chen-Ruan cohomology is for almost
complex orbifolds, in other words, the real dimensions of the
orbiolds are even. In this paper we study the similar Chen-Ruan
cohomology theory for almost contact orbifolds which have odd real
dimensions. The motivation for us to study this theory for almost
contact orbifolds is that the almost contact orbifolds are very
related to almost complex orbifolds.

An almost contact orbifold $X$ of dimension $2n+1$ is an orbifold
$X$ together with an almost contact structure. And the actions of
local groups of the orbifold preserve the almost contact structure.
The most interesting feature of the orbifold cohomology for the
almost contact orbifold is the cohomology of so-called twisted
sectors defined by [8],[2]. For a twisted sector $X_{(g)}$, the
property of almost contact structure implies that the action of
$g$ on the tangent space belongs to $U(n)\times 1$. Similar to the
definition of the degree shifting number in [2], we define the
degree shifting number $\iota_{(g)}$ of $X_{(g)}$. We have the
following definition:
\begin{defn}
Let $X$ be an almost contact orbifold of dimension $2n+1$, we
define the Chen-Ruan orbifold cohomology group of $X$ by
$$H_{orb}^{d}(X;\mathbb{Q}):=\bigoplus_{(g)\in
T_{1}}H^{d-2\iota_{(g)}}(X_{(g)},\mathbb{Q})$$
\end{defn}
where $T_{1}$ represents all the equivalent classes of $(g)$,
which show the number of twisted sectors. For the almost contact
orbifold $X$, $X\times \mathbb{R}$ is an almost complex orbifold.
We call it the corresponding almost complex orbiofld of $X$. Let
$H_{CR}^{d}(X\times \mathbb{R};\mathbb{Q})$ represents the
Chen-Ruan cohomology group defined in [2], we have the following
result:
\begin{thm}
Let $X$ be an almost contact orbifold of dimension $2n+1$, and
$X\times \mathbb{R}$ be its corresponding almost complex orbifold,
then as vector spaces, $H_{orb}^{d}(X;\mathbb{Q}) \cong
H_{CR}^{d}(X\times \mathbb{R};\mathbb{Q})$.
\end{thm}
To define the orbifold cup product on
$H_{orb}^{*}(X;\mathbb{Q})=\bigoplus_{d}H_{orb}^{d}(X;\mathbb{Q})$,
we still adopt the same method of Chen and Ruan [2] using the
degree zero Gromov-Witten invariant of the orbifold stable maps.
Let $\eta_{j}\in H^{d_{j}}(X_{(g_{j})};\mathbb{Q})$ for $j=1,2$,
$\eta_{3}\in H_{c}^{d_{3}}(X_{(g_{3})};\mathbb{Q})$ . we define
$<\eta_{1}\cup_{orb}\eta_{2},\eta_{3}>_{orb}:=
<\eta_{1},\eta_{2},\eta_{3}>_{orb}$, where
$$<\eta_{1},\eta_{2},\eta_{3}>_{orb}:=\sum_{{\bf g}\in T_{3}^{0}}\int_{X_{(\bf{g})}}^{orb}e^{*}_{1}\eta_{1}\wedge
 e^{*}_{2}\eta_{2}\wedge e^{*}_{3}\eta_{3}\wedge e_{A}(E_{(\bf{g})})
$$
$X_{({\bf g})}=X_{(g_{1},g_{2},g_{3})}$ is a 3-multisector with
$g_{1}g_{2}g_{3}=1$. The maps $e_{j}: X_{(\bf{g})}\longrightarrow
X_{(g_{j})}$ defined by $(p,(\mbox{\bf{g}})_{p})\longmapsto
(p,(g_{j})_{p})$ are the evaluation maps. $e_{A}(E_{(\bf{g})})$ is
the Euler form of the obstruction bundle $E_{({\bf g})}$
calculated from the connection $A$. The obstruction bundle
$E_{({\bf g})}$ over $X_{({\bf g})}$ is defined from the
obstruction bundle of the corresponding 3-multisector of its
corresponding almost complex orbiofld $X\times \mathbb{R}$. Under
this orbifold cup product, we have:
\begin{thm}
Let $X$ be an almost contact orbifold of dimension $2n+1$,  then
the orbifold cohomology $H_{orb}^{*}(X;\mathbb{Q})$ with the
orbifold cup product defined above is a cohomological ring.
\end{thm}

\begin{cor}
Let $X$ be an almost contact orbifold of dimension $2n+1$, and
$X\times \mathbb{R}$ be its corresponding almost complex orbifold,
then as cohomological rings, $H_{orb}^{*}(X;\mathbb{Q}) \cong
H_{CR}^{*}(X\times \mathbb{R};\mathbb{Q})$.
\end{cor}

The paper is organized  as follows. Section 2 is a review of some
basic facts concerning orbifold, orbifold vector bundle and almost
contact orbifold.  In section 3 we define the orbifold
cohomological ring of any almost contact orbifold. we introduce
twisted sectors, degree shifting numbers and orbifold cup product
for the almost contact orbifolds. Using the Chen-Ruan cup product
of the corresponding almost complex orbifolds of the almost
contact orbifolds, we prove that the orbifold cup product makes
the orbifold cohomology of an almost contact orbifold into a
cohomological ring. Finally in section 4 we give two examples.

\subsection*{Acknowledgments}
We would like to  thank Professor Yongbin Ruan  for very helpful
encouragements. And we especially thank him for discussing  this
interesting problem with us.


\section{Preliminaries}
\subsection{Orbifold and orbifold vector bundle.}
For the concept of orbifold, we know that the definition  of
Satake [11] is for reduced orbifolds.  Here we introduce the
general orbifolds which are not necessarily reduced. We give the
definition from Chen and Ruan [2].
\begin{defn}An orbifold structure on a
Hausdorff,  separate topological space $X$ is given by an open
cover $\mathcal{U}$ of $X$ satisfying the following conditions.

(1) Each element $U$ in $\mathcal{U}$ is uniformized, say by
$(V,G,\pi)$.  Namely, $V$ is a smooth manifold and $G$ is a finite
group acting smoothly on $V$ such that $U=V/G$ with $\pi$ as the
quotient map. Let $Ker(G)$ be the subgroup of $G$ acting trivially
on $V$.

(2) For $U^{'}\subset U$,  there is a collection of injections
$(V^{'},G^{'},\pi^{'})\longrightarrow (V,G,\pi)$. Namely, the
inclusion $i:U^{'}\subset U$ can be lifted to maps $\widetilde{i}:
V^{'}\longrightarrow V $ and an injective homomorphism $i_{*}:
G^{'}\longrightarrow G$ such that $i_{*}$ is an isomorphism from
$Ker(G^{'})$ to $Ker(G)$ and  $\widetilde{i}$ is
$i_{*}$-equivariant.

(3) For any point $x\in U_{1}\cap U_{2}$, $U_{1},U_{2}\in
\mathcal{U} $, there is a $U_{3}\in \mathcal{U} $ such that $x\in
 U_{3} \subset U_{1}\cap U_{2}$.
\end{defn}

For any point $x\in X$, suppose that $(V,G,\pi)$ is a uniformizing
neighborhood and $\overline{x}\in \pi^{-1}(x)$. Let $G_{x}$ be the
stabilizer of $G$ at $\overline{x}$. Up to conjugation, it is
independent of the choice of $\overline{x}$ and is called the
$local~ group$ of $x$. Then there exists a sufficiently small
neighborhood $V_{x}$ of $\overline{x}$ such that
$(V_{x},G_{x},\pi_{x})$ uniformizes a small neighborhood of $x$,
where $\pi_{x}$ is the restriction $\pi\mid V_{x}$.
$(V_{x},G_{x},\pi_{x})$ is called $a$ $local$ $ chart$ at $x$. The
orbifold structure is called $reduced $ if the action of $G_{x}$
is effective for every $x$.

Let $pr: E\longrightarrow X$ be a rank $k$ complex $orbifold~
bundle$ over an orbifold $X$([2]). Then a uniformizing  system for
$E\mid U=pr^{-1}(U)$ over a uniformized subset $U$ of $X$ consists
of the following data:

(1) A uniformizing system $(V,G,\pi)$ of $U$.

(2) A uniformizing system $(V \times \mathbb{C}^{k},$ $ G,
\widetilde{\pi})$ for $E\mid U$. The action of $G$ on $V\times
\mathbb{C}^{k}$ is an extension of the action of $G$ on $V$ given
by $g\cdot(x,v)=(g\cdot x,\rho(x,g)v)$ where $\rho: V\times
G\longrightarrow Aut(\mathbb{C}^{k})$ is a smooth map satisfying:
$$\rho(g\cdot x,h)\circ \rho(x,g)=\rho(x,hg),g,h\in G, x\in V.$$

(3)The natural projection map $\widetilde{pr}: V \times
\mathbb{C}^{k}\longrightarrow V$ satisfies $\pi\circ
\widetilde{pr}=pr\circ \widetilde{\pi}$.

By an orbifold connection $\bigtriangleup$ on $E$ we mean an
equivariant connection that satisfies
$\bigtriangleup=g^{-1}\bigtriangleup g$ for every uniformizing
system of $E$. Such a connection can be always obtained by
averaging an equivariant partition of unity.

A $C^{\infty}$ differential form on $X$ is a $G$-invariant
differential form on $V$ for each uniformizing system $(V,G,\pi)$.
Then orbifold integration is defined as follows. Suppose $U=V/G$
is connected, for any compactly supported differential $n$-form
$\omega$ on $U$, which is, by definition, a $G$-invariant $n$-form
$\widetilde{\omega}$ on $V$,

$$
 \int_{U}^{orb}\omega
:=\frac{1}{|G|}\int_{V}\widetilde{\omega}  \eqno{(2.1)}
$$

Where $|G|$ is the order of $G$.  The orbifold integration over
$X$ is defined by using a $C^{\infty}$ partition of unity. The
orbifold  integration coincides with the usual measure theoretic
integration iff the orbifold structure is reduced.


\subsection{The Almost Contact Orbifolds.}
In this section we introduce the basic concept of almost contact
orbifolds. First we give the definition of an almost contact
orbifold.
\begin{defn}
An almost  contact structure on a orbifold $X$ of dimension $2n+1$
is a co-dimension one distribution $\xi$ on $X$ which is locally
given by the kernel of a 1-form $\eta$ such that $\xi$ has an
almost complex structure $J$, and for every uniformizing system
$(V,G,\pi)$, the action of $G$ on $V$ preserves this almost
complex structure $J$.
\end{defn}

\begin{defn}
A contact structure on a orbifold $X$ of dimension $2n+1$ is a
co-dimension one distribution $\xi$ on $X$ which is given by the
kernel of a global 1-form $\eta$ with $\eta\wedge d\eta\neq 0$.
$\eta$ is said to represent $\xi$, and is called a contact form.
\end{defn}

For the next part of the paper, an almost contact orbifold is
always an orbifold such that the distribution $\xi$ has an almost
complex structure. So  the tangent bundle of $X$ has the structure
group $U(n)\times 1$. Given an almost contact orbifold $X$. For an
open subset $U$ of $X$, let $(V,G,\pi)$ be its local uniformizing
neighborhood. If $g\in G$ is an element of $G$, then the action of
$g$ on the open set $V$ belongs to $U(n)\times 1$.
\begin{rem}
For a contact orbifold $X$ with contact form $\eta$, the
symplectization of $X$ is a symplectic orbifold.  From [13], for
the co-oriented orbifold $X$, its symplectization can be
identified with $(X\times \mathbb{R},d(e^{t}\eta))$, where $t$ is
the $\mathbb{R}$ coordinate.
\end{rem}

\section{The Chen-Ruan Cohomology Rings of Almost Contact Orbifolds}
\subsection{Twisted Sectors and the Chen-Ruan Cohomology Group.}

First we introduce twisted sectors. Let $X$ be an orbifold.
Consider the set of pairs:
$$\widetilde{X}_{k}=\{(p,(\mbox{\bf{g}})_{G_{p}})|p\in X,
\mbox{\bf{g}}=(g_{1},\cdots,g_{k}),g_{i}\in G_{p}\}$$ where
$(\mbox{\bf{g}})_{G_{p}}$ is the conjugacy class of $k$-tuple
$\mbox{\bf{g}}=(g_{1},\cdots,g_{k})$ in $G_{p}$. We use $G^{k}$ to
denote the set of $k$-tuples. If there is no confusion, we will
omit the subscript $G_{p}$ to simplify the notation. Suppose that
$X$ has an orbifold structure $\mathcal{U}$ with uniformizing
systems $\{(\widetilde{U},G_{U},\pi_{U})\}$. From Chen and Ruan
[2], also see [8],   $\widetilde{X}_{k}$ is naturally an
orbifold,  with the generalized orbifold structure at
$(p,(\mbox{\bf{g}})_{G_{p}})$ given by
$(V_{p}^{\bf{g}},C(\mbox{\bf{g}}),\pi:V_{p}^{\bf{g}}\longrightarrow
V_{p}^{\bf{g}}/ C(\mbox{\bf{g}}))$,  where
$V_{p}^{\bf{g}}=V_{p}^{g_{1}}\cap \cdots V_{p}^{g_{k}}$,
$C(\mbox{\bf{g}})=C(g_{1})\cap\cdots C(g_{k})$. Here
$\mbox{\bf{g}}=(g_{1}, \cdots, g_{k})$,  $V_{p}^{g}$ stands for
the fixed point set of $g$ in $V_{p}$. When $X$ is a almost
contact orbifold of dimension $2n+1$, $\widetilde{X_{k}}$ inherits
a almost contact structure from $X$, and when $X$ is closed,
$\widetilde{X_{k}}$ is finite disjoint union of closed orbifolds.

Now we describe the the connected components of
$\widetilde{X_{k}}$,  Recall that every point $p$ has a local
chart $(V_{p},G_{p},\pi_{p})$ which gives a local uniformized
neighborhood $U_{p}=\pi_{p}(V_{p})$. If $q\in U_{p}$,  up to
conjugation there is a unique injective homomorphism $i_{*}:
G_{q}\longrightarrow G_{p}$.  For $\mbox{\bf{g}}\in (G_{q})^{k}$,
the conjugation class $i_{*}(\mbox{\bf{g}})_{q}$ is well defined.
We define an equivalence relation $i_{*}(\mbox{\bf{g}})_{q}\cong
(\mbox{\bf{g}})_{q}$. Let $T_{k}$ denote the set of equivalence
classes.To abuse the notation,  we use $(\mbox{\bf{g}})$ to denote
the equivalence class which $(\mbox{\bf{g}})_{q}$ belongs to.  We
will usually denote an element of $T_{1}$ by $(g)$. It is clear
that $\widetilde{X_{k}}$ can be decomposed as a disjoint union of
connected components:

$$\widetilde{X_{k}}=\bigsqcup_{({\bf g}) \in T_{k}}X_{(\bf{g})}$$

Where $X_{(\bf{g})}=\{(p,(\mbox{\bf{g}}^{'})_{p})|
\mbox{\bf{g}}^{'} \in (G_{p})^{k},(\mbox{\bf{g}}^{'})_{p}\in
(\mbox{\bf{g}})\}$. Note that for ${\bf g}=(1,\cdots,1)$, we have
$X_{(\bf{g})}=X$.  A component $X_{(\bf{g})}$ is called a
$k-multisector$, if $\bf{g}$ is not the identity. A component of
$X_{(g)}$ is simply called a $twisted$ $sector$. If $X$ has a
almost contact structure or an almost complex, then $X_{(\bf{g})}$
has the analogous structure induced from $X$. We define
$$T_{3}^{0}=\left\{({\bf g})=(g_{1},g_{2},g_{3})\in T_{3}| g_{1}g_{2}g_{3}=1\right\}.$$
Note that there is an  one to one correspondence between $T_{2}$
and $T_{3}^{0}$ given by $(g_{1},g_{2})\longmapsto
(g_{1},g_{2},(g_{1}g_{2})^{-1})$.

Now we define the Chen-Ruan cohomology group. Let $X$ be an almost
contact orbifold of dimension $2n+1$. Then for a point $p$ with
nontrivial group $G_{p}$, from section 2.2, we know that the
almost contact structure on $X$ gives rise to an effective
representation $\rho_{p}:G_{p}\longrightarrow
U(n,\mathbb{C})\times 1$. For any $g\in G_{p}$, we write
$\rho_{p}(g)$, up to conjugation, as a diagonal matrix
$$diag\left(e^{2\pi i\frac{m_{1,g}}{m_{g}}},\ldots,e^{2\pi i\frac{m_{n,g}}{m_{g}}}\right)\times 1.$$
where $m_{g}$ is the order of $g$ in $G_{p}$, and $0\leq m_{i,g}<
m_{g}$. Define a function $\iota: \widetilde{X_{1}}\longrightarrow
\mathbb{Q} $ by
$$
\iota(p,(g)_{p})=\sum_{i=1}^{n}\frac{m_{i,g}}{m_{g}}.
$$
We can see that the function $\iota:
\widetilde{X_{1}}\longrightarrow \mathbb{Q} $ is locally constant
and $\iota=0$ if $g=1$. Denote  its value on $X_{(g)}$ by
$\iota_{g}$. We call $\iota_{g}$ the degree shifting number of
$X_{(g)}$. It has the following properties:

(1) $\iota_{(g)}$ is an integer iff $\rho_{p}(g)\in
SL(n,\mathbb{C})\times 1$;

(2) $2(\iota_{(g)}+\iota_{(g^{-1})})=2rank(\rho_{p}(g)-Id\times 1)
=2n-dim_{\mathbb{R}}X_{(g)}$.

\begin{defn} Let $X$ be an almost contact orbifold of dimension $2n+1$,
we define the orbifold cohomology group of $X$ by

$$H_{orb}^{d}(X; \mathbb{Q}):=\bigoplus_{(g)\in
T_{1}}H^{d-2\iota_{(g)}}(X_{(g)},\mathbb{Q})$$

\end{defn}

Let $X\times \mathbb{R}$ be the corresponding almost complex
orbifold of $X$, then it is easy to see that the orbifold
structure of $X\times \mathbb{R}$ can be easily described. Let
$Y=X\times \mathbb{R}$, then a k-multisector $Y_{({\bf
g})}=X_{({\bf g})}\times \mathbb{R}$ for ${\bf
g}=(g_{1},\cdots,g_{k})$. We have the following theorem:
\begin{thm}
Let $X$ be an almost contact orbifold of dimension $2n+1$, and
$X\times \mathbb{R}$ be its corresponding almost complex orbifold
, then as vector spaces, $H_{orb}^{d}(X;\mathbb{Q}) \cong
H_{CR}^{d}(X\times \mathbb{R};\mathbb{Q})$.
\end{thm}
\begin{proof}
We see that there exists an one to one correspondence between the
twisted sectors of the orbifolds $X$ and $X\times \mathbb{R}$, and
this correspondence is given by $X_{({\bf g})}\longrightarrow
Y_{({\bf g})}=X_{({\bf g})}\times \mathbb{R}$. The degree shifting
numbers of the two twisted sectors are equal, and any dimensional
cohomology groups of $X$ and $X\times \mathbb{R}$ are isomorphic.
From the definition of the Chen-Ruan cohomology group of $X\times
\mathbb{R}$ in [2] and definition 3.1.1, the theorem is proved.
\end{proof}

\subsection{The Obstruction Bundle.}
Let $X$ be an almost contact orbifold of dimension $2n+1$, we
consider its corresponding almost complex orbifold  $Y=X\times
\mathbb{R}$. For a 3-mutisector $X_{({\bf
g})}=X_{(g_{1},g_{2},g_{3})}$ with $(g_{1},g_{2},g_{3})\in
T_{3}^{0}$, we know that $Y_{({\bf g})}=X_{({\bf g})}\times
\mathbb{R}$. Let $((p,x),({\bf g})_{(p,x)})$ be a generic point in
$X_{({\bf g})}\times \mathbb{R}$. Let $K({\bf g})$ be the subgroup
of $G_{p}$ generated by $g_{1}$ and $g_{2}$. Consider an orbifold
Riemann sphere with three orbifold points
$(S^{2},(p_{1},p_{2},p_{3}),(k_{1},k_{2},$ $k_{3}))$. When there
is no confusion, we will simply denote it by $S^{2}$. The orbifold
fundamental group is:

$$\pi_{1}^{orb}(S^{2})=\{\lambda_{1},\lambda_{2},\lambda_{3}|\lambda_{i}^{k_{i}}=1,\lambda_{1}\lambda_{2}\lambda_{3}=1\}$$

Where $\lambda_{i}$ is represented by a loop around the marked
$p_{i}$. There is a surjective homomorphism
$$\rho: \pi_{1}^{orb}(S^{2})\longrightarrow K(\bf{g})$$
specified by mapping $\lambda_{i}\longmapsto g_{i}$. $Ker(\rho)$
is a finite-index subgroup of $\pi_{1}^{orb}(S^{2})$. Let
$\widetilde{\Sigma}$ be the orbifold universal cover of $S^{2}$.
Let $\Sigma=\widetilde{\Sigma}/Ker(\rho)$. Then $\Sigma$ is
smooth,  compact and $\Sigma/K({\bf g})=S^{2}$. The genus of
$\Sigma$ can be computed using Riemann Hurwitz formula for Euler
characteristics of a branched covering,  and turns out to be

$$g(\Sigma)=\frac{1}{2}(2+|K(\mbox{\bf{g}})|-\Sigma_{i=1}^{3}\frac{|K(\bf{g})|}{k_{i}}) \eqno{(3.1)}$$

$K(\bf{g})$ acts holomorphically on $\Sigma$  and hence
$K(\bf{g})$ acts on $H^{0,1}(\Sigma)$. The "obstruction bundle"
$\widetilde{E}_{(\bf{g})}$ over $X_{({\bf g})}\times \mathbb{R}$
is constructed as follows. On the local chart
$(V_{p}^{\bf{g}}\times \mathbb{R},C(\bf{g}),\pi)$ of $X_{({\bf
g})}\times \mathbb{R}$, $\widetilde{E}_{(\bf{g})}$ is given by
$(T(V_{p}\times \mathbb{R})\otimes
H^{0,1}(\Sigma))^{K(\bf{g})}\times (V_{p}^{\bf{g}}\times
\mathbb{R}) \longrightarrow V_{p}^{\bf{g}}\times \mathbb{R}$,
where $(T(V_{p}\times \mathbb{R})\otimes
H^{0,1}(\Sigma))^{K(\bf{g})}$ is the $K(\bf{g})$-invariant
subspace. We define an action of $C(\bf{g})$ on $T(V_{p}\times
\mathbb{R})\otimes H^{0,1}(\Sigma)$, which is the usual one on
$T(V_{p}\times \mathbb{R})$ and trivial on $H^{0,1}(\Sigma)$. The
the action of $C(\bf{g})$ and $K(\bf{g})$ commute and
$(T(V_{p}\times \mathbb{R})\otimes H^{0,1}(\Sigma))^{K(\bf{g})}$
is invariant under $C(\bf{g})$. Thus we have obtained an action of
$C(\bf{g})$ on $(T(V_{p}\times \mathbb{R})\otimes
H^{0,1}(\Sigma))^{K(\bf{g})}\times (V_{p}^{\bf{g}}\times
\mathbb{R})\longrightarrow V_{p}^{\bf{g}}\times \mathbb{R}$,
extending the usual one on $V_{p}^{\bf{g}}\times \mathbb{R}$.
These trivializations fit together to define the bundle
$\widetilde{E}_{(\bf{g})}$ over $X_{(\bf{g})}\times \mathbb{R}$.
If we set $e\times 1 : X_{({\bf g})}\times
\mathbb{R}\longrightarrow X\times \mathbb{R}$ to be the map given
by $((p,x),({\bf g})_{(p,x)})\longmapsto (p,x)$, one may think of
$\widetilde{E}_{(\bf{g})}$ as $((e\times 1)^{*}T(X\times
\mathbb{R})\otimes H^{0,1}(\Sigma))^{K(\bf{g})}$. The rank of
$\widetilde{E}_{(\bf{g})}$  is given by the formula [2]:

$$rank_{\mathbb{R}}(\widetilde{E}_{(\bf{g})})=
dim_{\mathbb{R}}(X_{({\bf g})}\times
\mathbb{R})-dim_{\mathbb{R}}(X\times \mathbb{R}
)+2\Sigma_{j=1}^{3}\iota_{(g_{j})} \eqno{(3.2)}$$

Let $e: X_{(\bf{g})}\longrightarrow X$ to be the map given by
$(p,({\bf g})_{p})\longmapsto p$, we have the following
proposition:
\begin{prop}
Let $\pi_{({\bf g})}: X_{({\bf g})}\times
\mathbb{R}\longrightarrow X_{({\bf g})}$  be the natural
projection, then $\left((e\times 1)^{*}T(X\times
\mathbb{R})\otimes H^{0,1}(\Sigma)\right)^{K(\bf{g})} \cong
\pi_{({\bf g})}^{*}\left((e^{*}TX\otimes
H^{0,1}(\Sigma))^{K(\bf{g})}\right)$
\end{prop}
\begin{proof}
We have 
 \begin{eqnarray} ((e\times
1)^{*}T(X\times\mathbb{R})\otimes
H^{0,1}(\Sigma))^{K(\bf{g})}&=&\pi_{({\bf g})}^{*}\left((e^{*}TX\otimes
H^{0,1}(\Sigma))\oplus (T\mathbb{R}\otimes
H^{0,1}(\Sigma))\right)^{K({\bf g})} \nonumber \\
&=&\pi_{({\bf g})}^{*}\left((e^{*}TX\otimes H^{0,1}(\Sigma))^{K({\bf g})}\oplus
(T\mathbb{R}\otimes
H^{0,1}(\Sigma)\right)^{K({\bf g})}\nonumber \\
&=&\pi_{({\bf g})}^{*}\left(e^{*}TX\otimes H^{0,1}(\Sigma)\right)^{K({\bf g})}
\nonumber
\end{eqnarray}
The third equality is right because $K({\bf g})$ acts on
$T\mathbb{R}$ trivially, and $(H^{0,1}(\Sigma))^{K({\bf
g})}=H^{0,1}(S^{2})=0$. So we complete the proof of the
proposition.
\end{proof}

We define the obstruction bundle $E_{({\bf g})}$ over $X_{({\bf g}
)}$ as $\left(e^{*}TX\otimes H^{0,1}(\Sigma)\right)^{K({\bf g})}$.
It is easy to see that the rank of $E_{({\bf g})}$ is equal to the
rank of $\widetilde{E}_{({\bf g})}$.  The rank is:
$$rank_{\mathbb{R}}(E_{(\bf{g})})=
dim_{\mathbb{R}}(X_{({\bf
g})})-dim_{\mathbb{R}}(X)+2\Sigma_{j=1}^{3}\iota_{(g_{j})}
\eqno{(3.3)}$$

\subsection{The Chen-Ruan Cup Product.}

For an almost contact orbifold $X$,  there is a natural map $I:
X_{(g)}\longrightarrow X_{(g^{-1})}$ defined by
$(p,(g)_{p})\longmapsto (p,(g^{-1})_{p})$.

\begin{defn}(Poincar$\acute{e}$ Duality) Let $X$ be an almost contact orbifold of dimension $2n+1$.
For any  $0\leq d\leq 2n$, the pairing

$$<,>_{orb}: H^{d}_{orb}(X)\times H^{2n-d}_{orb,c}(X)\longrightarrow \mathbb{Q}$$
is defined by taking the direct sum of

$$<,>_{orb}^{(g)}: H^{d-2\iota_{(g)}}(X_{(g)};
\mathbb{Q})\times
H_{c}^{2n-d-2\iota_{(g^{-1})}}(X_{(g^{-1})};\mathbb{Q})
\longrightarrow \mathbb{Q}$$ where
$$<\alpha,\beta>_{orb}^{(g)}=\int_{X_{(g)}}^{orb}\alpha\wedge
I^{*}(\beta)$$ for $\alpha\in H^{d-2\iota_{(g)}}(X_{(g)};
\mathbb{Q})$, and $\beta\in
H^{2n-d-2\iota_{(g^{-1})}}(X_{(g^{-1})}; \mathbb{Q})$.
\end{defn}

Choose an orbifold connection $A$ on $E_{(\bf{g})}$. Let
$e_{A}(E_{(\bf{g})})$ be the Euler form computed from the
connection $A$ by Chen-Weil theory. Let $\eta_{j}\in
H^{d_{j}}(X_{(g_{j})}; \mathbb{Q})$, for $j=1,2$, $\eta_{3}\in
H_{c}^{d_{3}}(X_{(g_{3})}; \mathbb{Q})$. Define maps $e_{j}:
X_{(\bf{g})}\longrightarrow X_{(g_{j})}$ by $(p,({\bf
g})_{p})\longmapsto (p,(g_{j})_{p})$.

\begin{defn} Similar to the definition in [2], we
define the 3-point function to be
$$<\eta_{1},\eta_{2},\eta_{3}>_{orb}:=\sum_{{\bf g}\in T_{3}^{0}}\int_{X_{(\bf{g})}}^{orb}e^{*}_{1}\eta_{1}\wedge
 e^{*}_{2}\eta_{2}\wedge e^{*}_{3}\eta_{3}\wedge e_{A}(E_{(\bf{g})})
 \eqno{(3.4)}$$
\end{defn}

Note that the above integral does not depend on the choice of $A$.
As in the definition 3.1.1, we extend the 3-point function to
$H_{orb}^{*}(X)$ by linearity. We define the orbifold cup product
by the relation

$$<\eta_{1}\cup_{orb}\eta_{2},\eta_{3}>_{orb}:=
<\eta_{1},\eta_{2},\eta_{3}>_{orb}  \eqno{(3.5)}$$

Again we extend $\cup_{orb}$ to $H_{orb}^{*}(X)$ via linearity.
Note that if $({\bf g})=(1,1,1)$, then
$\eta_{1}\cup_{orb}\eta_{2}$ is just the ordinary cup product
$\eta_{1}\cup \eta_{2}$ in $H^{*}(X)$.

\subsection{Associativity.}
In this section we prove the associativity of the Chen-Ruan cup
product defined in the section 3.3. What we do here is to compare
the CR-cup product of $H^{*}_{orb}(X)$ with the orbifold cup
product of the chen-Ruan cohomology of its corresponding almost
complex orbifold $X\times \mathbb{R}$ defined in [2]. We have the
following theorem:
\begin{thm}
Let $X$ be an almost contact orbifold of dimension $2n+1$,  then
the orbifold cohomology $H_{orb}^{*}(X;\mathbb{Q})$ with the
Chen-Ruan cup product defined in section 3.3 is a cohomological
ring.
\end{thm}
\begin{proof}
From section 3.1, we see that there is an one to one
correspondence between the twisted sectors of $X$ and its
corresponding almost complex orbifold  $Y=X\times \mathbb{R}$. By
definition 3.1.1 and the definition of the orbifold cohomology of
$X\times \mathbb{R}$ defined in [2], we have:
$$H_{orb}^{*}(X;\mathbb{Q}):=\bigoplus_{(g)\in
T_{1}}H^{*-2\iota_{(g)}}(X_{(g)},\mathbb{Q})$$ and
$$H_{orb}^{*}(X\times \mathbb{R};\mathbb{Q}):=\bigoplus_{(g)\in
T_{1}}H^{*-2\iota_{(g)}}(X_{(g)}\times \mathbb{R},\mathbb{Q})$$
For the k-multisector $X_{({\bf g})}$, we know that $Y_{({\bf
g})}=X_{({\bf g})}\times \mathbb{R}$. Let $\pi_{({\bf g})}:
X_{({\bf g})}\times \mathbb{R}\longrightarrow X_{({\bf g})}$ be
the natural projection. Then from [1], $(\pi_{({\bf g})})_{*}:
H_{cv}^{*}(X_{({\bf g})}\times \mathbb{R})\longrightarrow
H^{*-1}(X_{({\bf g})})$ is an isomorphism.

Let $X_{({\bf g})}$ be a 3-multisector with $g_{1}g_{2}g_{3}=1$.
For any $\eta_{1}\in H^{*}(X_{(g_{1})})$, $\eta_{2}\in
H^{*}(X_{(g_{2})})$, $\eta_{3}=(\pi_{(g_{3})})_{*}\gamma\in
H^{*}(X_{(g_{3})})$, where $\gamma\in H_{cv}^{*}(X_{({\bf
g})}\times \mathbb{R})$. We have the following commutative
diagram:
\begin{center}
\setlength{\unitlength}{0.8cm}
\begin{picture}(5,3)\put(0,2){\shortstack{$X_{({\bf g})}\times \mathbb{R}$}}
\put(2.6,2.3){\shortstack{$\pi_{({\bf g})}$}}
 \put(2.0,2.1){\vector(1,0){2.0}}
 \put(4.2,2){\shortstack{$X_{({\bf g})}$}}
 \put(1.2,1.3){\shortstack{$\downarrow$}}
 \put(0,1.3){\shortstack{$e_{i}\times 1$}}
 \put(4.1,1.4){\shortstack{$e_{i}$}}
 \put(4.5,1.3){\shortstack{$\downarrow$}}
 \put(0 ,0.7){\shortstack{$X_{(g_{i})}\times \mathbb{R}$}}
\put(4.2,0.7){\shortstack{$X_{(g_{i})}$}}
 \put(2.0,0.8){\vector(1,0){2.0}}
 \put(2.6,0.99){\shortstack{$\pi_{(g_{i})}$}}
\end{picture}
\end{center}
for $i=1,2,3$. So we have $(e_{i}\times 1)^{*}(\pi_{(g_{i})})^{*}
\eta_{i}=(\pi_{({\bf g})})^{*}e_{i}^{*}\eta_{i}$ for i=1,2 and
$(\pi_{({\bf g})})_{*}(e_{3}\times
1)^{*}\gamma=e_{3}^{*}(\pi_{(g_{3})})_{*}\gamma$. From the
definition of the 3-point function for almost complex orbifolds in
[2] and the orbifold version of the proposition 6.15 in [1], we
have:
$$<(\pi_{(g_{1})})^{*}\eta_{1},(\pi_{(g_{2})})^{*}\eta_{2},
\gamma>_{orb}$$
\begin{eqnarray}
&=&\sum_{{\bf g}\in T_{3}^{0}}\int_{X_{(\bf{g})}\times \mathbb{R}
}^{orb}(e_{1}\times 1)^{*}(\pi_{(g_{1})})^{*}\eta_{1}\wedge
 (e_{2}\times 1)^{*}(\pi_{(g_{2})})^{*}\eta_{2}\wedge
 (e_{3}\times 1)^{*}\gamma\wedge
 e_{\widetilde{A}}(\widetilde{E}_{(\bf{g})})\nonumber \\
 &=&\sum_{{\bf g}\in T_{3}^{0}}\int_{X_{(\bf{g})}\times \mathbb{R}
}^{orb}(\pi_{({\bf g})})^{*}e_{1}^{*}\eta_{1}\wedge
 (\pi_{({\bf g})})^{*}e_{2}^{*}\eta_{2}\wedge
 (e_{3}\times 1)^{*}\gamma\wedge
 (\pi_{({\bf g})})^{*}e_{A}(E_{(\bf{g})})\nonumber \\
&=&(-1)^{deg\gamma\cdot deg e_{A}(E_{{\bf g}})}\sum_{{\bf g}\in T_{3}^{0}}\int_{X_{(\bf{g})}\times \mathbb{R}
}^{orb}(\pi_{({\bf g})})^{*}(e_{1}^{*}\eta_{1}\wedge
 e_{2}^{*}\eta_{2}\wedge
e_{A}(E_{(\bf{g})}))\wedge (e_{3}\times 1)^{*}\gamma \nonumber \\
&=&(-1)^{deg\gamma\cdot deg e_{A}(E_{{\bf g}})}\sum_{{\bf g}\in T_{3}^{0}}\int_{X_{(\bf{g})}}^{orb}
e_{1}^{*}\eta_{1}\wedge
 e_{2}^{*}\eta_{2}\wedge
e_{A}(E_{(\bf{g})})\wedge (\pi_{({\bf g})})_{*}(e_{3}\times 1)^{*}\gamma \nonumber \\
&=&-\sum_{{\bf g}\in T_{3}^{0}}\int_{X_{(\bf{g})}}^{orb}
e_{1}^{*}\eta_{1}\wedge
 e_{2}^{*}\eta_{2}\wedge
 e_{3}^{*}(\pi_{(g_{3})})_{*}\gamma\wedge e_{A}(E_{(\bf{g})}) \nonumber \\
&=&-<\eta_{1},\eta_{2},(\pi_{(g_{3})})_{*}\gamma>_{orb} \nonumber
\end{eqnarray}
We see that the orbifold cup product of the Chen-Ruan cohomology
$H^{*}_{orb}(X; \mathbb{Q})$ is identified with the CR-product of
the orbifold cohomology for the corresponding almost complex
orbifold $X\times \mathbb{R}$ of $X$ module a sign. From [2], the CR-product
for the almost complex orbifold $X\times \mathbb{R}$ satisfies the
associativity, so the orbifold cup product of the Chen-Ruan
cohomology $H^{*}_{orb}(X; \mathbb{Q})$ for the almost contact
orbifold $X$ also satisfies the associativity.  So the Chen-Ruan
cohomology $H^{*}_{orb}(X; \mathbb{Q})$ for the almost contact
orbifold $X$ with this orbifold cup product forms a cohomological
ring.  This complete the proof of the theorem.
\end{proof}
From the proof of the above theorem, we have the following
corollary:
\begin{cor}
Let $X$ be an almost contact orbifold of dimension $2n+1$, and
$X\times \mathbb{R}$ be its corresponding almost complex orbifold,
then as cohomological rings, $H_{orb}^{*}(X;\mathbb{Q}) \cong
H_{CR}^{*}(X\times \mathbb{R};\mathbb{Q})$.
\end{cor}

\section{Examples.}
\subsection{Example 4.1}
The 3-sphere $X=S^{3}$ is an almost contact manifold with standard
almost contact form.  Let $\mathbb{Z}_{3}$ be a cyclic group of
order $3$. suppose that $\mathbb{Z}_{3}$ acts on $S^{3}$ in a
plane and preserves the almost contact form, then the quotient
$S^{3}/\mathbb{Z}_{3}$ is a almost contact orbifold of dimension
$3$. We see that this orbifold has one singular submanifold
$S^{1}$, and the local orbifold group in this singular set is the
cyclic group $\mathbb{Z}_{3}$. Let $g=e^{2\pi i \frac{1}{3}}\in
\mathbb{Z}_{3}$, then $X_{(g)}=S^{1}$ and   $X_{(g^{2})}=S^{1}$,
and the degree shifting numbers for the two twisted sectors are
$\frac{1}{3}$ and $\frac{2}{3}$. So from the definition 3.1.1, we
have:
$$H^{d}_{orb}(S^{3}/\mathbb{Z}_{3};\mathbb{Q})=
H^{d}(S^{3}/\mathbb{Z}_{3};\mathbb{Q})\bigoplus
H^{d-\frac{2}{3}}(S^{1};\mathbb{Q}) \bigoplus
H^{d-\frac{4}{3}}(S^{1};\mathbb{Q})$$ For ${\bf g}_{1}=(g,g,g)\in
T_{3}^{0}$, we have $X_{({\bf g})}=S^{1}$, and $X_{{\bf
g}_{2}}=X_{(g^{2},g^{2},g^{2})}=S^{1}$. From the formula (3.3), we
see that the rank of the obstruction bundle $E_{({\bf g}_{1} )}$
over $X_{({\bf g}_{1})}$ is 0, the rank of the obstruction bundle
$E_{({\bf g}_{2} )}$ over $X_{({\bf g}_{2})}$ is 2. So for
$\eta_{i}\in H^{*}(X_{(g)};\mathbb{Q}) (i=1,2,3)$, from section
3.3, we have
$$<\eta_{1},\eta_{2},\eta_{3}>_{orb}:=\int_{X_{({\bf g}_{1})}}^{orb}e^{*}_{1}\eta_{1}\wedge
 e^{*}_{2}\eta_{2}\wedge e^{*}_{3}\eta_{3}
$$
The integration is the usual integration on the orbifold and we
can calculate easily.

For $\eta_{i}\in H^{*}(X_{(g^{2})};\mathbb{Q}) (i=1,2,3)$, from
section 3.3, we have
$$<\eta_{1},\eta_{2},\eta_{3}>_{orb}:=\int_{X_{({\bf g}_{2})}}^{orb}e^{*}_{1}\eta_{1}\wedge
 e^{*}_{2}\eta_{2}\wedge e^{*}_{3}\eta_{3}\wedge e_{A}(E_{({\bf g}_{2})})
$$
But the dimension of the Euler class $e_{A}(E_{({\bf g}_{2})})$ is
$2$, while $X_{({\bf g}_{2})}$ has dimension $1$, so the above
3-point function is zero.

For the 3-multisector $X_{(g,g^{2},1)}$, it is easy to calculate
that the obstruction bundle over it has dimension zero. The
3-point function is the usual integration on orbifold and we can
easily calculate it. So we complete the analysis the Chen-Ruan
cohomology ring of the almost contact orbifold
$S^{3}/\mathbb{Z}_{3}$.

\subsection{Example 4.2}
We consider the weighted projective space
$X=\mathbb{P}(1,2,2,3,3,3)$. From the example of [7], we know
that $X$ is an orbifold. And the twisted sectors of this orbifold
are  $X_{(g_{2})}=\mathbb{P}(2,2)$,
$X_{(g_{1})}=\mathbb{P}(3,3,3)$ and
$X_{(g_{1}^{2})}=\mathbb{P}(3,3,3)$, the degree shifting numbers
are $\iota_{(g_{1})}=\frac{5}{3}$,
$\iota_{(g_{1}^{2})}=\frac{4}{3}$ and $\iota_{(g_{2})}=2$. Let
$S^{1}$ be the 1-dimensional sphere. Then $Y=X\times
S^{1}=\mathbb{P}(1,2,2,3,3,3)\times S^{1}$ is an almost contact
orbifold. It is easy to see that the twisted sectors of $Y$ are
$Y_{(g_{2}\times 1)}=\mathbb{P}(2,2)\times S^{1}$,
$Y_{(g_{1}\times 1)}=\mathbb{P}(3,3,3)\times S^{1}$ and
$Y_{(g_{1}^{2}\times 1)}=\mathbb{P}(3,3,3)\times S^{1}$ and the
degree shifting numbers are $\iota_{(g_{1}\times 1)}=\frac{5}{3}$,
$\iota_{(g_{1}^{2}\times 1)}=\frac{4}{3}$ and $\iota_{(g_{2}\times
1 )}=2$. So the Chen-Ruan cohomology group of
$Y=\mathbb{P}(1,2,2,3,3,3)\times S^{1}$ is
\begin{eqnarray}
H^{d}_{orb}(\mathbb{P}(1,2,2,3,3,3)\times
S^{1};\mathbb{Q})\nonumber &=&
H^{d}(\mathbb{P}(1,2,2,3,3,3)\times S^{1};\mathbb{Q})\\
\nonumber &\bigoplus & H^{d-\frac{10}{3}}(\mathbb{P}(3,3,3)\times
S^{1};\mathbb{Q}) \\ \nonumber &\bigoplus &
H^{d-\frac{8}{3}}(\mathbb{P}(3,3,3)\times S^{1};\mathbb{Q})\\
\nonumber  &\bigoplus & H^{d-4}(\mathbb{P}(2,2)\times
S^{1};\mathbb{Q}) \nonumber
\end{eqnarray}
For the Chen-Ruan ring structure, we see that the 3-multisectors
are $Y_{(g_{1}\times 1,g_{1}\times 1,g_{1}\times 1)}$
$=Y_{(g_{1}^{2}\times 1,g_{1}^{2}\times 1,g_{1}^{2}\times
1)}=\mathbb{P}(3,3,3)\times S^{1}$, $Y_{(g_{1}\times
1,g_{1}^{2}\times 1,1)}=\mathbb{P}(3,3,3)\times S^{1}$ and
$Y_{(g_{2}\times 1,g_{2}\times 1,1)}$ $=\mathbb{P}(2,2)\times
S^{1}$. For the 3-multisectors $Y_{(g_{1}\times 1,g_{1}^{2}\times
1,1)}$ and $Y_{(g_{2}\times 1,g_{2}\times 1,1)}$, it is easily to
see that the dimensions of the obstruction bundles over these
3-multisectors are zero, so the integrations of the 3-point
function (3.4) are the usual integration on the orbifolds, we can
determine them easily.

For the 3-multisector $Y_{(g_{1}\times 1,g_{1}\times 1,g_{1}\times
1)}$ $=\mathbb{P}(3,3,3)\times S^{1}$, let $\eta_{i}\in
H^{*}(\mathbb{P}(3,3,3)\times S^{1};\mathbb{Q})$ for $i=1,2,3$.
From (3.3), the dimension of the obstruction bundle $E_{({\bf
g}\times 1)}$ over $Y_{(g_{1}\times 1,g_{1}\times 1,g_{1}\times
1)}$ is 4. The integration
$$<\eta_{1},\eta_{2},\eta_{3}>_{orb}:=\int_{\mathbb{P}(3,3,3)\times S^{1}}^{orb}e^{*}_{1}\eta_{1}\wedge
 e^{*}_{2}\eta_{2}\wedge e^{*}_{3}\eta_{3}\wedge e_{A}(E_{({\bf
g}\times 1)}) \eqno{(4.1)}$$ is nonzero only if there is  some
$\eta_{j}\in H^{1}(\mathbb{P}(3,3,3)\times S^{1};\mathbb{Q})$. We
assume that $\eta_{1}\in H^{1}(\mathbb{P}(3,3,3)\times
S^{1};\mathbb{Q})$,$\eta_{2}\in H^{0}(\mathbb{P}(3,3,3)\times
S^{1};\mathbb{Q})$, $\eta_{3}\in H^{0}(\mathbb{P}(3,3,3)\times
S^{1};\mathbb{Q})$. So the integration $(4.1)$ is
$$<\eta_{1},\eta_{2},\eta_{3}>_{orb}:=\int_{\mathbb{P}(3,3,3)\times S^{1}}^{orb}e^{*}_{1}\eta_{1}
\wedge e_{A}(E_{({\bf g}\times 1)}) \eqno{(4.2)}$$ Because
$H^{1}(\mathbb{P}(3,3,3)\times
S^{1};\mathbb{Q})=H^{0}(\mathbb{P}(3,3,3);\mathbb{Q})\otimes
H^{1}(S^{1};\mathbb{Q})$, $e_{1}^{*}\eta_{1}$ represents the class
of $H^{1}(S^{1};\mathbb{Q})$. $e_{1}^{*}(\eta_{1})$ is the
poncar$\acute{e}$ duality of the suborbifold $S^{1}$ in
$\mathbb{P}(3,3,3)\otimes S^{1}$, $e_{A}(E_{({\bf g}\times
1)})=e_{A}(E_{({\bf g})})$, so from [1], we have:
$$\int_{\mathbb{P}(3,3,3)\times S^{1}}^{orb}e^{*}_{1}\eta_{1}
\wedge e_{A}(E_{({\bf g}\times 1)})=\int_{\mathbb{P}(3,3,3)}^{orb}
e_{A}(E_{({\bf g})})  \eqno{(4.3)}$$ The right side in (4.3) is
calculated in the example of [7]. So the formula (4.2) can be
calculated.

For the 3-multisector $Y_{(g_{1}^{2}\times 1,g_{1}^{2}\times
1,g_{1}^{2}\times 1)}$ $=\mathbb{P}(3,3,3)\times S^{1}$, let
$\eta_{i}\in H^{*}(\mathbb{P}(3,3,3)\times S^{1};\mathbb{Q})$ for
$i=1,2,3$. From (3.3), the dimension of the obstruction bundle
$E_{({\bf g}\times 1)}$ over $Y_{(g_{1}^{2}\times
1,g_{1}^{2}\times 1,g_{1}^{2}\times 1)}$ is 2. The integration
$$<\eta_{1},\eta_{2},\eta_{3}>_{orb}:=\int_{\mathbb{P}(3,3,3)\times S^{1}}^{orb}e^{*}_{1}\eta_{1}\wedge
 e^{*}_{2}\eta_{2}\wedge e^{*}_{3}\eta_{3}\wedge e_{A}(E_{({\bf
g}\times 1)}) \eqno{(4.4)}$$ is nonzero only if there are  some
$\eta_{i}\in H^{1}(\mathbb{P}(3,3,3)\times S^{1};\mathbb{Q})$,
$\eta_{j}\in H^{2}(\mathbb{P}(3,3,3)\times S^{1};\mathbb{Q})$,
because $H^{1}(\mathbb{P}(3,3,3)\times
S^{1};\mathbb{Q})=H^{0}(\mathbb{P}(3,3,3);\mathbb{Q})\otimes
H^{1}(S^{1};\mathbb{Q})$ and $H^{2}(\mathbb{P}(3,3,3)\times
S^{1};\mathbb{Q})=H^{2}(\mathbb{P}(3,3,3);\mathbb{Q})\otimes
H^{0}(S^{1};\mathbb{Q})$. Let $\eta_{1}\in
H^{1}(\mathbb{P}(3,3,3)\times S^{1};\mathbb{Q})$,$\eta_{2}\in
H^{2}(\mathbb{P}(3,3,3)\times S^{1};\mathbb{Q})$, $\eta_{3}\in
H^{0}(\mathbb{P}(3,3,3)\times S^{1};\mathbb{Q})$. So the
integration $(4.4)$ is
$$<\eta_{1},\eta_{2},\eta_{3}>_{orb}:=\int_{\mathbb{P}(3,3,3)\times S^{1}}^{orb}e^{*}_{1}\eta_{1}
\wedge e_{2}^{*}\eta_{2}\wedge e_{A}(E_{({\bf g}\times 1)})
\eqno{(4.5)}$$ We know that $e_{1}^{*}\eta_{1}$ represents the
class of $H^{1}(S^{1};\mathbb{Q})$. $e_{1}^{*}(\eta_{1})$ is the
poncar$\acute{e}$ duality of the suborbifold $S^{1}$ in
$\mathbb{P}(3,3,3)\otimes S^{1}$, $e_{A}(E_{({\bf g}\times
1)})=e_{A}(E_{({\bf g})})$, so from [1], we have:
$$\int_{\mathbb{P}(3,3,3)\times S^{1}}^{orb}e^{*}_{1}\eta_{1}
\wedge e_{2}^{*}\eta_{2}\wedge e_{A}(E_{({\bf g}\times
1)})=\int_{\mathbb{P}(3,3,3)}^{orb} e_{2}^{*}\eta_{2}\wedge
e_{A}(E_{({\bf g})}) \eqno{(4.6)}$$ The right side in (4.6) is
calculated in the example of [7]. So the formula (4.5) can be
calculated. We complete the analysis of the Chen-Ruan ring
structure of the orbifold $\mathbb{P}(1,2,2,3,3,3)\times S^{1}$.


\subsection*{}

\bibliographystyle{amsplain}
\bibliography{xibi}

\end{document}